\documentclass[11pt]{article}
\usepackage{amsmath,amssymb}
\usepackage[latin1]{inputenc}
\usepackage{epsfig}
\numberwithin{equation}{section}
\newtheorem{thm}{Theorem}[section]

\newtheorem{aprop}[thm]{Proposition}
\newtheorem{acor}[thm]{Corollary}
\newtheorem{arem}[thm]{Remark}
\newenvironment{adem}[1][]%
   {\ \\ {\bf Proof #1~: }}%
   {\hfill\mbox{\rule{2 true mm}{3 true mm}}\vskip 2 ex\noindent}

\newcommand{\E}{{\mathbb E}}
\newcommand{\R}{{\mathbb R}}
\renewcommand{\P}{{\mathbb P}}

\title{Stochastic flows approach to Dupire's formula}
 \author{B.Jourdain\thanks{CERMICS, project-team Mathfi, \'Ecole des
     Ponts, ParisTech, 6-8 av Blaise Pascal, Cit\'e
     Descartes, Champs sur Marne, 77455 Marne-la-Vall\'ee Cedex 2, France -
     e-mail : jourdain@cermics.enpc.fr }}

\begin{document}
 \maketitle  
\begin{abstract}The probabilistic equivalent formulation of Dupire's
  PDE is the Put-Call duality equality. In local volatility models including exponential L\'evy jumps, we give a direct
  probabilistic proof for this result based on stochastic flows
  arguments. This approach also enables us to check the probabilistic
  equivalent formulation of various generalizations of Dupire's PDE recently
  obtained by Pironneau \cite{pironneau} by the adjoint equation
  technique in the case of complex options. 
\end{abstract}
\section*{Introduction}
The second order derivative of the price of a Call option with respect
to the strike variable is equal to the risk-neutral density of the underlying
stock at maturity multiplied by the actualization factor.
In a stock model with a local volatility function and a proportional dividend
rate (\eqref{edsvoloc} with $\mu=m=0$), Dupire \cite{dup} takes advantage
of this specificity to obtain a PDE (see \eqref{edpdup} for $m=0$) satisfied by the Call pricing
function in the maturity and strike variables. His proof consists in
integrating twice in space the 
Fokker-Planck equation governing the time evolution of the density of
the stock price. Alternatively, one may use the
Green function of the problem or the
adjoint equation technique \cite{pironneau}. For calibration purposes,
Dupire's PDE permits to express the local volatility function in terms
of the function giving the Call prices for all strikes and maturities.\par
Dupire's PDE can be interpretated as the pricing PDE for a Put
option. This leads to the Put-Call duality (equality \eqref{dual} for
$\tilde{\mu}=\mu=m=0$) which is in fact an equivalent
formulation : the Call price is
transformed into the Put price by simultaneous exchange of the interest
and dividend rates and of the spot and strike prices in addition to
time-reversal of the local volatility function. To our knowledge, no
direct probabilistic proof is available for the equality of the expectations
giving the Call and Put prices. In \cite{andreasencarr}, in models
including exponential L\'evy jumps, Carr and Andreasen derive a PIDE generalizing Dupire's PDE by computing the evolution of the
Call payoff with respect to maturity thanks to the It\^o-Tanaka formula and
taking expectations. The present paper deals with such models (see
\eqref{edsvoloc}). In the second section, we
give a probabilistic proof of the Put-Call duality \eqref{dual} equality equivalent
to this PIDE. We check the equality of the expectations by an
argument based on stochastic flows of diffeomorphisms. The flow
properties of the SDE \eqref{edsvoloc} involved in this argument are
introduced in the first section and proved in the appendix.\par
In a recent paper, Pironneau \cite{pironneau} obtains various
generalizations of Dupire's PDE to complex options by the adjoint equation
  technique. More precisely, for a given complex option, he shows that
  it is possible to compute the pricing function for all strikes and maturities by solving a single PDE. In calibration procedures, solving
  this PDE instead of one pricing PDE for the maturity and strike of
  each quoted option permits important computation time reduction. Most
  of these generalized Dupire's PDEs have equivalent probabilistic
  interpretations similar to the Put-Call duality. In the third and
  fourth sections of the paper, we use stochastic flows to check the equivalent interpretations
  corresponding to binary and options written on two assets.\par The fifth section deals
  with barrier options in local volatility models without jumps. In section 1.1 \cite{pironneau}, Pironneau addresses two-barriers options but we have only been able to give a probabilistic
  equivalent interpretation (see \eqref{dualbar}) in the one-barrier case. Moreover,
  besides particular cases, it seems challenging to give a probabilistic
  proof of this equivalent formulation. The case of American options is not addressed in \cite{pironneau}. In
\cite{alfonjour}, we deal with the case of perpetual options when the
local volatility function does not depend on time. For the
perpetual American Call price to be equal to the perpetual American Put price, in
addition to the exchanges of the interest and dividend rates and of the
spot and strike prices, the volatility function has to be modified. Our
approach consists in deriving and studying an ODE satisfied by the exercise boundary
as a function of the strike variable. Again, a direct probabilistic proof
of the duality results appears challenging. The stochastic flow approach
presented in the present paper does not seem suited to deal with options
involving stopping times like barrier and American options.

{\bf Notations :} For $T>0$ and $m$ a measure on $\R$ such that $\int_\R (1+e^l)\wedge
l^2m(dl)<+\infty$,
let $(W_t)_{t\in[0,T]}$ be a standard Brownian motion and $\mu$ denote an
independent Poisson random point measure on
$(0,T]\times \R$ with intensity $m(dl)dt$.\\We consider the following risk-neutral evolution for the underlying
stock price 
\begin{equation}
  dX^x_t=\sigma(t,X^x_t)X^x_tdW_t+(r-\delta)X^x_tdt+X^x_{t^-}\int_\R
   (e^l-1)(\mu(dt,dl)-m(dl)dt),\;X^x_0=x>0\label{edsvoloc} 
\end{equation}
where $r$ denotes the interest rate and $\delta$ the dividend rate. The
local volatility function $\sigma(t,x)$ is assumed to belong to the
space 
$${\mathcal V}=\left\{f:[0,T]\times(0,+\infty)\rightarrow
  \R:\;\sup_{[0,T]\times\R}\sum_{k=0}^3|x^k\partial^k_2f(t,x)|<+\infty\right\}$$
where $\partial_2^k f$ denotes the $k$-th order derivative of the
function $f$ with respect to its second variable.\\
The process
$(\hat{W}_t=W_{T-t}-W_T)_{t\in[0,T]}$ obtained by time-reversal of $W$
is a Brownian motion independent from the image $\hat{\mu}$ of $\mu$ by the
mapping $(t,l)\in(0,T]\times\R \rightarrow (T-t,-l)$ which is a Poisson
random point measure on $[0,T)\times\R$ with intensity $\hat{m}(dl)dt$ where
$\hat{m}$ denotes the image of $m$ by $l\in\R\rightarrow -l$.\\Let us
also introduce the L\'evy processes
\begin{align}
   &L_t=\int_{(0,t]\times\R} l(\mu(ds,dl)-1_{\{|l|\leq 1\}}m(dl)ds)-t\int_\R (e^l-1-l1_{\{|l|\leq 1\}})m(dl)\notag\\&\mbox{and
}\hat{L}_t=L_{T-t}-L_T=\int_{[0,t)\times\R} l(\hat{\mu}(ds,dl)-1_{\{|l|\leq
  1\}}\hat{m}(dl)ds)+t\int_\R (e^l-1-l1_{\{|l|\leq 1\}})m(dl).\label{levy}
\end{align}
By the L\'evy-Kinchine formula, 
\begin{align}
   \E(e^{iuL_t})=e^{t\psi(u)}\mbox{ where }\psi(u)&=-iu\int_\R(e^l-1-l1_{\{|l|\leq
  1\}})m(dl)+\int_\R(e^{iul}-1-iul1_{\{|l|\leq
  1\}})m(dl)\notag\\&=\int_\R(e^{iul}-1-iu(e^l-1))m(dl).\label{levykin}
\end{align}

\section{Stochastic flows of diffeomorphisms}
According to the theory of stochastic flows of diffeomorphisms
developped by Kunita \cite{kunita}, almost surely, the solution at time
$t>0$ of a Stochastic Differential Equation with regular coefficients is
a diffeomorphism as a function of the initial position. Derivatives of
the solution with respect to the initial condition solve the linear
equations obtained by formal derivation of the SDE. Last, the inverse
diffeomorphism is associated with a backward SDE. In the following
proposition, we adapt these results to a slight generalization of the
SDE with jumps preserving positivity \eqref{edsvoloc} considered in the present paper.
\begin{aprop}\label{propflot2}
Assume that $\sigma,\beta\in{\mathcal V}$ and let $\eta(t,x)=x\sigma(t,x)$. 
   Then trajectorial uniqueness holds for the stochastic
   differential equations
\begin{align*}
   &dX^x_t=\eta(t,X^x_t)dW_t+\beta(t,X^x_t)X^x_tdt+X^x_{t^-}\int_\R
   (e^l-1)(\mu(dt,dl)-m(dl)dt),\;t\leq T,\;X^x_0=x>0\\
&dZ^z_t=\eta(T-t,Z^z_t)d\hat{W}_t+Z^z_t(\sigma\partial_2\eta-\beta)(T-t,Z^z_t)dt+Z^z_{t^-}\int_\R(e^l-1)(\hat{\mu}(dt,dl)+m(dl)dt),\;t\leq
T\end{align*}
where $Z^z_0=z>0$ and $\int_\R(e^l-1)(\hat{\mu}(dt,dl)+m(dl)dt)$ stands
  for $$\int_\R(e^l-1)(\hat{\mu}(dt,dl)-1_{|l|\leq
  1}\hat{m}(dl)dt)+\int_\R(e^l-1+1_{\{|l|\leq 1\}}(e^{-l}-1))m(dl).$$
They admit solutions such that for almost all $\omega\in\Omega$, 
the mappings $x\rightarrow X^x_T$ and $z\rightarrow Z^z_T$ are inverse
increasing diffeomorphisms of $(0,+\infty)$,
\begin{align}&\forall
t\in[0,T],\;Z^z_t=ze^{\int_0^t\sigma(T-s,Z^z_s)d\hat{W}_s+\int_0^t[\sigma\partial_2\eta-\beta-\frac{\sigma^2}{2}](T-s,Z^z_s)ds+\hat{L}_{t^+}}\mbox{
  (with $\hat{L}_{T^+}=\hat{L}_T$)}\label{soledsback},\\
   &\partial_x
   X^x_T=e^{\int_0^T\partial_2\eta(s,X^x_s)dW_s+\int_0^T(\beta+X^x_s\partial_2\beta-\frac{(\partial_2\eta)^2}{2})(s,X^x_s)ds+L_T}\label{edsder}\\&\mbox{and}\;\forall x,z>0,\;\{X^x_T\geq z\}=\{x\geq Z^z_T\}.\label{ensemblegaux}
\end{align}
\end{aprop}
The rather technical proof of this proposition is postponed to the
appendix. To deduce the Put-Call duality equality \eqref{dual}, we are
going to check the equality of the derivatives of both sides with
respect to $x$. The next result enables us to justify the formula
$\partial_x
\E(e^{-rT}(X^x_T-y)^+)=\E\left(e^{-rT}\partial_xX^x_T1_{\{X^x_T\geq
    y\}}\right)$ obtained by formal derivation and where the indicator
funtion in the right-hand side will be replaced thanks to
\eqref{ensemblegaux}. Its proof is also postponed to the appendix.
\begin{aprop} Under the assumptions and notations of Proposition
  \ref{propflot2}, when for some $z>0$, the local volatility function
  $\sigma$ does not vanish on a neighbourhood of $(T,z)$ in $[0,T]\times
  (0,+\infty)$, then
\begin{equation}
   \forall x>0,\;\P(X^x_T=z)=\P(Z^z_T=x)=0.\label{points}
\end{equation}
Last, if $\beta(t,x)=\gamma$ for some constant $\gamma\in\R$ then 
\begin{equation}
   \forall x>0,\;\E(e^{-\gamma T}X^x_T)=x\mbox{ and }\E(e^{-\gamma
   T}\partial_x X^x_T)=1
\label{moment}\end{equation}
and for any sequence $(h_n)_{n\geq 0}$ of non-zero real numbers greater
   than $-x$ converging to
   zero, the random variables
   $\left((X^{x+h_n}_T-X^x_T)/h_n\right)_{n\geq 0}$ are uniformly
   integrable.
\label{propflot3}\end{aprop}

\section{Standard options}
For $y>0$, let $C(T,x,y)={\mathbb E}(e^{-rT}(X^x_T-y)^+)$ denote the
price of the Call option with
maturity $T$ and strike $y$ written on the underlying $X^x$ evolving
according to \eqref{edsvoloc}.\\
We are going to derive the Put-Call duality \eqref{dual} from the
following proposition.
\begin{thm}\label{dupder}
Assume that the local volatility function does not vanish in a
neighbourhood of $(T,y)$ in $[0,T]\times (0,+\infty)$. Then
   \begin{equation}
      \forall x>0,\;\partial_x C(T,x,y)=\partial_x{\mathbb
  E}\left(e^{-\delta T}(x-Y^y_T)^+\right)\label{dualder}
   \end{equation}
where
$dY^y_t=\sigma(T-t,Y^y_t)Y^y_td\hat{W}_t+(\delta-r)Y^y_tdt+Y^y_{t^-}\int_\R(e^{-l}-1)(\tilde{\mu}(dt,dl)-e^lm(dl)dt),\;Y^y_0=y$
and $\tilde{\mu}$ is a Poisson random point measure with intensity
$e^lm(dl)dt$ independent from $\hat{W}$.
\end{thm}
As ${\mathbb
  E}\left(e^{-\delta T}(x-Y^y_T)^+\right)\leq e^{-\delta
  T}x$ and $C(T,x,y)\leq {\mathbb E}(e^{-rT}X^x_T)$ with ${\mathbb
  E}(e^{-rT}X^x_T)=e^{-\delta T}x$ according to \eqref{moment}, one has $\lim_{x\rightarrow 0^+}C(T,x,y)=\lim_{x\rightarrow 0^+}{\mathbb
  E}\left(e^{-\delta T}(x-Y^y_T)^+\right)=0$. Hence the Put-Call duality
  follows from \eqref{dualder} :
\begin{acor}If the local volatility function does not vanish in a
neighbourhood of $(T,y)$\label{dupus},
\begin{equation}
   \forall x>0,\;C(T,x,y)={\mathbb
  E}\left(e^{-\delta T}(x-Y^y_T)^+\right).\label{dual}
\end{equation}
\end{acor}
\begin{arem}\label{equivdupdual}If the local volatility
  function $\sigma$ is positive and belongs
  to ${\mathcal V}$ for any $T>0$, then the Put-Call
  duality \eqref{dual} holds for all
  $(T,y)\in[0,+\infty)\times(0,+\infty)$. This implies Dupire's PIDE. Indeed, for
$s\in[0,T]$, let $(Y^{s,y,T}_t)_{t\in[s,T]}$ solve $Y^{s,y,T}_s=y$ and for
  $t\in[s,T]$,
$$dY^{s,y,T}_t=\sigma(T-t,Y^{s,y,T}_t)Y^{s,y,T}_td\hat{W}_t+(\delta-r)Y^{s,y,T}_tdt+Y^{s,y,T}_{t^-}\int_\R(e^{-l}-1)(\tilde{\mu}(dt,dl)-e^lm(dl)dt).$$
For fixed $s$, by time translation, the expectation
  $\E\left(e^{-\delta(T-(T-s))}(x-Y^{T-s,y,T}_T)^+\right)$ does
  not depend on $T$ greater than $s$ and may be denoted $u(s,y)$. By the Feynman-Kac formula, one has
\begin{equation*}
   \begin{cases}
      \partial_s u(s,y)=\frac{1}{2}\sigma^2(s,y)y^2\partial_{yy}
      u(s,y)+(\delta-r)y\partial_yu(s,y)-\delta
      u(s,y)\\\phantom{\partial_s u(s,y)=}+\int_\R \left[u(s,ye^{-l})-u(s,y)-\partial_y u(s,y)y(e^{-l}-1)\right]e^lm(dl),\;s,y>0\\
u(0,y)=(x-y)^+,\;y>0
   \end{cases}.
\end{equation*}
Since $C(T,x,y)=u(T,y)$, one deduces that the function $C$ solves
Dupire's PIDE in the variables $(T,y)$ :
   \begin{equation}
   \begin{cases}
      \partial_T C(T,x,y)=\frac{1}{2}\sigma^2(T,y)y^2\partial_{yy}
      C(T,x,y)+(\delta-r)y\partial_yC(T,x,y)-\delta C(T,x,y)\\\phantom{\partial_T C(T,x,y)=}+\int_\R \left[C(T,x,ye^{-l})-C(T,x,y)-\partial_y C(T,x,y)y(e^{-l}-1)\right]e^lm(dl),\;T,y>0\\
C(0,x,y)=(x-y)^+,\;y>0
   \end{cases}.\label{edpdup}
\end{equation}
Conversely, if \eqref{edpdup} holds, the function $v(t,y)=C(T-t,x,y)$
   satisfies the PIDE
\begin{align*}
   \partial_t v(t,y)&+\frac{1}{2}\sigma^2(T-t,y)y^2\partial_{yy}
      v+(\delta-r)y\partial_yv-\delta v\\&+\int_\R \left[v(t,ye^{-l})-v(t,y)-\partial_y v(t,y)y(e^{-l}-1)\right]e^lm(dl)=0
\end{align*}
with terminal condition $v(T,y)=(y-x)^+$. By the Feynman-Kac
      representation for the solution of this PIDE, $v(0,y)=\E\left(e^{-\delta
      T}(x-Y^y_T)^+\right)$ and \eqref{dual} holds.
\end{arem}
\begin{adem}[of Theorem \ref{dupder}]
Let $(h_n)_{n\geq 0}$ be a sequence of non-zero real numbers greater
   than $-x$ converging to
   zero. Since $x\rightarrow X^x_T$ is increasing according to
   Proposition \ref{propflot2}, one has
$$0\leq \frac{(X^{x+h_n}_T-y)^+-(X^x_T-y)^+}{h_n}\leq
\frac{X^{x+h_n}_T-X^x_T}{h_n},$$
and the variables
$\left(((X^{x+h_n}_T-y)^+-(X^x_T-y)^+)/h_n\right)_{n\geq 0}$ are
uniformly integrable by Proposition \ref{propflot3}. By \eqref{points}, these variables converge a.s. to $\partial_x
X^x_T1_{\{X^x_T\geq y\}}$ as $n\rightarrow +\infty$. One deduces that
$C(T,x,y)$ is differentiable with respect to $x$ and
$\partial_xC(T,x,y)=e^{-\delta T}\E\left(e^{(\delta-r)T}\partial_x
  X^x_T1_{\{X^x_T\geq y\}}\right)$. By \eqref{edsder}, this implies
$$\partial_xC(T,x,y)=e^{-\delta
  T}\E\left(e^{L_T}\E\left(e^{\int_0^T\partial_2\eta(u,X^x_u)dW_u-\frac{1}{2}\int_0^T(\partial_2\eta(u,X^x_u))^2du}1_{\{X^x_T\geq y\}}\bigg|(L_s)_{s\in[0,T]}\right)\right).$$
Since $\partial_2\eta(t,x)=\sigma(t,x)+x\partial_2\sigma(t,x)$ is
bounded on $[0,T]\times\R$, for
$y=0$ the conditional expectation in the
right-hand-side is equal to $1$ by Novikov's criterion (Proposition
1.15 p.308 \cite{reyor}). Therefore, by Girsanov theorem,
$$\E\left(e^{\int_0^T\partial_2\eta(u,X^x_u)dW_u-\frac{1}{2}\int_0^T(\partial_2\eta(u,X^x_u))^2du}1_{\{X^x_T\geq y\}}\bigg|(L_s)_{s\in[0,T]}\right)=\E\left(1_{\{\bar{X}^x_T\geq
      y\}}\bigg|(L_s)_{s\in[0,T]}\right)$$
$$\mbox{where
}d\bar{X}^x_t=\eta(t,\bar{X}^x_t)dW_t+\beta(t,\bar{X}^x_t)\bar{X}^x_tdt+\bar{X}^x_{t^-}\int_\R
l(\mu(dl,dt)-1_{\{|l|\leq
  1\}}m(dl)dt),\;\bar{X}^x_0=x$$
with $\beta(t,z)=\sigma\partial_2\eta(t,z)+(r-\delta)$. By
\eqref{ensemblegaux}, one deduces that
\begin{equation}
   \partial_x C(T,x,y)=e^{-\delta T}\E\left(e^{L_T}1_{\{x\geq Z^y_T\}}\right),\label{exprder}
 \end{equation} where,
according to \eqref{soledsback} and the definition of $\beta$,
$$\forall
t\in[0,T],\;Z^z_t=ze^{\int_0^t\sigma(T-s,Z^z_s)d\hat{W}_s+(\delta-r)t-\frac{1}{2}\int_0^t\sigma^2(T-s,Z^z_s)ds+\hat{L}_{t^+}}\mbox{
    (convention : $\hat{L}_{T^+}=\hat{L}_T$)}.$$ 
When $m=0$, we are done. Otherwise, we still have to derive the dynamics of
$\hat{L}_{t^+}$ under the probability measure with density $e^{L_T}$
with respect to $\P$.
By \eqref{levykin}, for $t\in[0,T]$ and $u\in\R$,
$\E(e^{L_T}e^{iu\hat{L}_{t^+}})=\E(e^{(1-iu)(L_T-L_{T-t})})\E(e^{L_{T-t}})=e^{t\psi(-(u+i))}$
and 
\begin{align*}
  \psi(-(u+i))&=\int_\R(e^le^{-iul}-1+(iu-1)(e^l-1))m(dl)\\&=\int_\R(e^{-iul}-1+iul1_{\{|l|\leq 1\}})e^lm(dl)+iu\int_\R(e^l-1-le^l1_{\{|l|\leq 1\}})m(dl).
\end{align*}
Therefore, under the probability measure with density $e^{L_T}$ with
respect to $\P$, 
$$\hat{L}_{t^+}=-\int_{(0,t]\times\R}l(\bar{\mu}(ds,dl)-1_{\{|l|\leq
  1\}}e^lm(dl)ds)+t\int_\R(e^l-1-le^l1_{\{|l|\leq 1\}})m(dl)$$
with $\bar{\mu}$ a Poisson random point measure with intensity
$e^lm(dl)dt$ independent from $W$. By It\^o's formula
$$dZ^z_t=\eta(T-t,Z^z_t)d\hat{W}_t+(\delta-r)Z^z_tdt+Z^z_{t^-}\int_\R
(e^{-l}-1)(\bar{\mu}(dt,dl)-e^lm(dl)dt),\;Z^z_0=z.$$
Since trajectorial and therefore weak uniqueness holds for this SDE with
jumps, $\E(e^{L_T}1_{\{x\geq Z^y_T\}})=\P(x\geq Y^y_T)$ and by \eqref{exprder},
$\partial_x
C(T,x,y)=e^{-\delta T}\E(1_{\{x\geq Y^y_T\}})$. 
According to
Lebesgue's Theorem, the
right-hand-side is equal to $\partial_x\E\left(e^{-\delta
    T}(x-Y^y_T)^+\right)$ since $\P(Y^y_T=x)=\E(e^{L_T}1_{\{Z^y_T=x\}})=0$ by \eqref{points}. 
\end{adem}
\begin{arem}\label{girs}Many authors have obtained another type of
  Put-Call duality by the following change of num\'eraire approach :
$$\E(e^{-rT}(X^x_T-y)^+)=\E^{\mathbb Q}\left(e^{-\delta
    T}\left(x-\frac{yx}{X^x_T}\right)^+\right)\mbox{ with
}\frac{d{\mathbb Q}}{d\P}=e^{(\delta-r)T}\frac{X^x_T}{x}.$$
In exponential L\'evy models, the local volatility function is constant
and Fajardo and Mordecki \cite{fajmor} derive \eqref{dual} by this approach. But in general, the coefficients of the SDE with jumps satisfied by
$\frac{xy}{X^x_t}$ under $\mathbb Q$ depend on the primal spot variable
and dual strike variable $x$. Then it is not clear to deduce a PIDE in
the variables $T$ and $y$ for the pricing function $C(T,x,y)$.
\end{arem}
\section{Binary options}
For $z>0$, let $c(T,x,z)=\E\left(e^{-rT}1_{\{X^x_T\geq
    z\}}\right)$ denote the price of the binary Call option with strike $z$ and maturity
    $T$ written on the underlying $X^x$. The following result is a
    direct consequence of \eqref{ensemblegaux} :
\begin{aprop}\label{dupbin}
   \begin{equation}
   \forall x,z>0,\;c(T,x,z)=\E\left(e^{-rT}1_{\{x\geq
    Z^z_T\}}\right)\label{dualbin}
\end{equation} 
where $Z^z$ solves 
  $dZ^z_t=\eta(T-t,Z^z_t)d\hat{W}_t+(\eta\partial_2\eta(T-t,Z^z_t)+(\delta-r)Z^z_t)dt+Z^z_{t^-}\int_\R(e^l-1)(\hat{\mu}(dt,dl)+m(dl)dt),\;Z^z_0=z.$\end{aprop}
\begin{arem}\label{pdebin}
   By a reasoning similar to the one made for the
  standard Dupire's PIDE in Remark \ref{equivdupdual}, one deduces that
  as a function of the maturity and strike variables, the function $c(T,x,z)$ solves
\begin{equation}
   \begin{cases}
    \partial_T
    c(T,x,z)=\frac{1}{2}\sigma^2(T,z)z^2\partial_{zz}c(T,x,z)+(\eta\partial_2\eta(T,z)+(\delta-r)z)\partial_zc(T,x,z)-rc(T,x,z)\\\phantom{\partial_T c(T,x,z)=}+\int_\R \left[c(T,x,ze^{-l})-c(T,x,z)+\partial_zc(T,x,z)z(e^l-1)\right]m(dl),\;T,z>0\\
c(0,x,z)=1_{\{x\geq z\}},\;z>0  
   \end{cases}.\label{edpbin}
\end{equation}
Remarking that $$\frac{1}{2}\sigma^2(T,z)z^2\partial_{zz}
      c(T,x,z)+\eta\partial_2\eta(T,z)\partial_{z}
      c(T,x,z)=\frac{1}{2}\partial_z\left(\sigma^2(T,z)z^2\partial_zc(T,x,z)\right),$$
      one recognizes the PDE obtained in \cite{pironneau} in the absence
      of jumps ($m$=0).
\end{arem}
\begin{arem}
   The standard Call and the binary Call pricing functions are linked by
$$C(T,x,y)=\E\left(e^{-rT}\int_y^{+\infty}1_{\{X^x_T\geq
    z\}}dz\right)=\int_y^{+\infty}c(T,x,z)dz.$$
Integrating the PIDE \eqref{edpbin} with respect to $z$ on the
    interval $[y,+\infty[$ and remarking that 
\begin{align*}
   \int_y^{+\infty}&\left[c(T,x,ze^{-l})-c(T,x,z)+\partial_zc(T,x,z)z(e^l-1)\right]dz\\&=\int_{ye^{-l}}^{+\infty}c(T,x,w)e^{l}dw+(e^l-1)\int_y^{+\infty}\partial_z(zc(T,x,z))dz-e^l\int_y^{+\infty}c(T,x,z)dz\\
&=-\left[e^lC(T,x,ye^{-l})+(e^l-1)y\partial_yC(T,x,y)-e^lC(T,x,y)\right]\\&=-e^l\left[C(T,x,ye^{-l})-C(T,x,y)-(e^{-l}-1)y\partial_yC(T,x,y)\right]
\end{align*}one recovers \eqref{edpdup}. This
    alternative proof of \eqref{edpdup} relies on properties of the
    derivative of the pricing function $C$ with respect to the strike variable,
    whereas Proposition \ref{dupder} deals with its derivative with
    respect to the spot variable.\end{arem}
\section{Options written on two assets}
We now consider a model with two assets evolving according to the
following dynamics
\begin{equation}
   \begin{cases}
      dX^{1,x_1}_t=\sigma_1(t,X^{1,x_1}_t)X^{1,x_1}_tdW^1_t+(r-\delta_1)X^{1,x_1}_tdt\\
\phantom{dX^{1,x_1}_t=}+X^{1,x_1}_{t^-}\int_{\R^2}(e^{l_1}-1)(\mu(dt,dl_1,\R)-m(dl_1,\R)dt),\;X^{1,x_1}_0=x_1\\
dX^{2,x_2}_t=\sigma_2(t,X^{2,x_2}_t)X^{2,x_2}_tdW^2_t+(r-\delta_2)X^{2,x_2}_tdt\\\phantom{dX^{2,x_2}_t=}+X^{2,x_2}_{t^-}\int_{\R^2}(e^{l_2}-1)(\mu(dt,\R,dl_2)-m(\R,dl_2)dt),\;X^{2,x_2}_0=x_2
   \end{cases},\label{mod2d}
\end{equation}
where for $i\in\{1,2\}$, $\sigma_i\in{\mathcal V}$. Here $W^1$ and $W^2$ are two standard Brownian motions such that
$d<W^1,W^2>_t=\rho_t dt$ with $\rho$ an adapted process and $\mu$ is an
independent Poisson random point measure on $(0,T]\times\R^2$ with
intensity $m(dl_1,dl_2)dt$ where
$$\int_{\R^2}(e^{l_1}+e^{l_2}+1)\wedge(l_1^2+l_2^2)m(dl_1,dl_2)<+\infty.$$
Let for $i\in\{1,2\}$, $\eta_i(t,z_i)=z_i\sigma_i(t,z_i)$,
  $(\hat{W}^i_t=W^i_{T-t}-W^i_T)_{t\in[0,T]}$ and $\hat{\mu}$ denote the
  image of $\mu$ by the mapping
  $(t,l_1,l_2)\in(0,T]\times\R^2\rightarrow (T-t,-l_1,-l_2)$.
\begin{aprop}Let $(Z^{1,z_1},Z^{2,z_2})$ solve
\begin{align*}
   dZ^{1,z_1}_t=\eta_1(T-t,Z^{1,z_1})d\hat{W}^1_t&+(\eta_1\partial_2\eta_1(T-t,Z^{1,z_1}_t)+(\delta_1-r)Z^{1,z_1}_t)dt\\
&+Z^{1,z_1}_{t^-}\int_{\R^2}(e^{l_1}-1)(\hat{\mu}(dt,dl_1,\R)+m(dl_1,\R)dt),\;Z^{1,z_1}_0=z_1\\
dZ^{2,z_2}_t=\eta_2(T-t,Z^{2,z_2})d\hat{W}^2_t&+(\eta_2\partial_2\eta_2(T-t,Z^{2,z_2}_t)+(\delta_2-r)Z^{2,z_2}_t)dt\\
&+Z^{2,z_2}_{t^-}\int_{\R^2}(e^{l_2}-1)(\hat{\mu}(dt,\R,dl_2)+m(\R,dl_2)dt),\;Z^{2,z_2}_0=z_2.
\end{align*}
Then for $w(T,x_1,x_2,z_1,z_2)=\E\left(e^{-rT}1_{\{x_1\geq Z^{1,z_1}_T
    ,x_2\geq Z^{2,z_2}_T\}}\right),$ one has 
\begin{align}
   &\forall
    x_1,x_2,y>0,\E\left((X^{1,x_1}_T+X^{2,x_2}_T-y)^+\right)=\int_0^yw(T,x_1,x_2,z,y-z)dz\notag\\&\phantom{\forall
    T,x_1,x_2,y>0,\;\E\left((X^{1,x_1}_T+X^{2,x_2}_T-y)^+\right)}+\int_y^{+\infty}w(T,x_1,x_2,z,0)+w(T,x_1,x_2,0,z)dz\label{dualbask}\\
&\mbox{and }\E\left((X^{1,x_1}_T\vee X^{2,x_2}_T-y)^+\right)=\int_y^{+\infty}w(T,x_1,x_2,z,0)+w(T,x_1,x_2,0,z)-w(T,x_1,x_2,z,z)dz.\notag
\end{align}
\end{aprop}
\begin{arem}
   When $\rho_t=-q(t,X^{1,x_1}_t,X^{2,x_2}_t)$, by a reasoning similar to
the one made in Remark \ref{pdebin}, one
checks that the function $w$ solves the following PIDE obtained in
\cite{pironneau} Section 2.1 in the absence of jumps :
\begin{equation*}
   \begin{cases}
    \partial_T
    w=\sum_{i=1}^2\left[\frac{1}{2}\partial_{z_i}\left(\sigma_i^2(T,z_i)z_i^2\partial_{z_i}w\right)+(\delta_i-r)z_i\partial_{z_i}w\right]-2qz_1z_2\partial_{z_1z_2}w-rw\\
\phantom{\partial_T
    w=}+\int_{\R^2}w(T,z_1e^{-l_1},z_2e^{-l_2})-[w-\sum_{i=1}^2z_i(e^{l_i}-1)\partial_{z_i}w](T,z_1,z_2)m(dl_1,dl_2)\\
w(0,x_1,x_2,z_1,z_2)=1_{\{x_1\geq z_1,x_2\geq z_2\}} 
   \end{cases}.\end{equation*}
\end{arem}
\begin{adem}
One has $$\forall y_1,y_2,y\geq 0,\;(y_1+y_2-y)^+=\int_0^y1_{\{y_1\geq z,y_2\geq
  y-z\}}dz+\int_y^{+\infty}1_{\{y_1\geq z\}}+1_{\{y_2\geq z\}}dz.$$
Therefore 
$$\E\left((X^{1,x_1}_T+X^{2,x_2}_T-y)^+\right)=e^{-rT}\E\left(\int_0^y1_{\{X^{1,x_1}_T\geq
    z,X^{2,x_2}_T\geq y-z\}}dz+\int_y^{+\infty}1_{\{X^{1,x_1}_T\geq
    z\}}+1_{\{X^{2,x_2}_T\geq
    z\}}dz\right).$$
Since by \eqref{ensemblegaux}, a.s., $\forall x_1,x_2,z_1,z_2>0,\;\{X^{1,x_1}_T\geq
    z_1\}=\{x_1\geq Z^{1,z_1}_T\}$ and $\{X^{2,x_2}_T\geq
    z_2\}=\{x_2\geq Z^{2,z_2}_T\}$,
    \begin{equation*}\;C(T,x_1,x_2,y)=e^{-rT}\E\left(\int_0^y1_{\{x_1\geq
    Z^{1,z}_T,x_2\geq Z^{2,y-z}_T\}}dz+\int_y^{+\infty}1_{\{x_1\geq
    Z^{1,z}_T\}}+1_{\{x_2\geq Z^{2,z}_T\}}dz\right)
 \end{equation*}
As by Proposition \ref{propflot2}, $z_2\rightarrow Z^{2,z_2}_T$
(resp. $z_1\rightarrow Z^{1,z_1}_T$) is an increasing
diffeomorphism of $(0,+\infty)$, $\lim_{z_2\rightarrow
  0^+}w(T,x_1,x_2,z_1,z_2)=e^{-rT}\P(x_1\geq Z^{1,z_1}_T)$ (resp. $\lim_{z_1\rightarrow
  0^+}w(T,x_1,x_2,z_1,z_2)=e^{-rT}\P(x_2\geq Z^{2,z_2}_T)$). One easily deduces \eqref{dualbask}.\\The
formula for the best-off Call option is obtained similarly remarking that
$$\forall y_1,y_2,y\geq 0,\;(y_1\vee y_2-y)^+=\int_y^{+\infty}1_{\{y_1\geq z\}}+1_{\{y_2\geq z\}}-1_{\{y_1\geq z,y_2\geq
  z\}}dz.$$
\end{adem}
\begin{arem}
   When $\sigma_i$ may depend smoothly on $x_j$ for $j\neq i$ in \eqref{mod2d}, $(x_1,x_2)\rightarrow
   (X^{1,x_1,x_2}_T,X^{2,x_1,x_2}_T)$ is still a diffeomorphism with
   inverse $(z_1,z_2)\rightarrow (Z^{1,z_1,z_2}_T,Z^{2,z_1,z_2}_T)$
   obtained as the solution of a two-dimensional SDE at time $T$. But in
   general, the events $\{X^{1,x_1,x_2}_T\geq
    z_1,X^{2,x_1,x_2}_T\geq
    z_2\}$ and $\{x_1\geq Z^{1,z_1,z_2}_T,x_2\geq Z^{2,z_1,z_2}_T\}$ are
   not equal and the argument above fails.
\end{arem}
\section{Barrier options}
 In absence of jumps $\boxed{m=0}$ and in the
 particular case of equal interest and dividend rates
 $\boxed{r=\delta}$, as a consequence of Proposition 1.2 \cite{pironneau}, 
\begin{equation}\label{dualbar}
   \forall x,y\geq
   z>0,\;\E\left((X^x_T-y)^+1_{\{\tau^x_z>T\}}\right)=\E\left((x-Y^y_{T\wedge t^y_z})^+\right)
\end{equation}
with $X^x$ solving \eqref{edsvoloc}, $Y^y$ defined in Proposition
\ref{dupder}, $\tau^x_z=\inf\{t\geq 0:X^x_t\leq z\}$ and
$t^y_z=\inf\{t\geq 0:Y^y_t\leq z\}$. Notice that since
$r=\delta$, there is no drift term in the dynamics of $X^x$ and
$Y^y$ and both processes are martingales in their natural
filtration.
This equality generalizes \eqref{dual} which can be recovered by taking
the limit $z\rightarrow 0$. It is easy to prove when either $x$ or $y$
is equal to $z$. Indeed, when $y\geq x=z$, both sides are equal to
$0$. And when $x\geq y=z$, using the martingale property of $X^x$ for
the third equality,
\begin{align*}
 \E\left((X^x_T-y)^+1_{\{\tau^x_y>T\}}\right)&=\E\left((X^x_T-y)1_{\{\tau^x_y>T\}}\right)=\E(X^x_T)-\E\left(X^x_T1_{\{\tau^x_y\leq
 T\}}\right)-y\P(\tau^x_y>T)\\&=x-y\P(\tau^x_y\leq
 T)-y\P(\tau^x_y>T)
\end{align*}
while the right-hand-side of \eqref{dualbar} is obviously equal to $x-y$.\par
Equation \eqref{dualbar} is equivalent to an equality where $x$ and $y$
play symmetric roles. Indeed substracting \eqref{dualbar} to
\eqref{dual}, and using the martingale property of $Y^y$ for the
second equality, one gets
\begin{align*}
 \E&\left((X^x_T-y)^+1_{\{\tau^x_z\leq T\}}\right)=\E\left(\left((x-Y^y_T)^+-(x-z)\right)1_{\{t^y_z\leq T\}}\right)\\&=\E\left(\left((x-Y^y_T)^+-(x-Y^y_T)\right)1_{\{t^y_z\leq T\}}\right)=\E\left((Y^y_T-x)^+1_{\{t^y_z\leq T\}}\right).
\end{align*}
In case the local volatility function $\sigma$ does not depend on the
time variable, this last equality obviously holds when $x=y$.\\
Beside these particular cases, it seems challenging to give a probabilistic proof of
\eqref{dualbar}. Derivation of the equality with respect to $x$ or $y$
does not lead to nice probabilistic equalities like \eqref{dualder}
obtained in the case of standard options.
Even in the Black-Scholes model with constant
volatility, \eqref{dualbar} is not obvious. For instance, the change of
num\'eraire approach presented in Remark \ref{girs} then enables to check that
$\E\left((X^x_T-y)^+1_{\{\tau^x_y>T\}}\right)=\E\left((x-Y^y_T)^+1_{\{\sup_{[0,T]}Y^y_t<\frac{yx}{z}\}}\right)$.
But it is not clear that the right-hand-side coincides with the one in
\eqref{dualbar}.
\section*{Appendix}
The proof of Proposition \ref{propflot2} relies on the following result
concerning SDEs without jumps which is a consequence of Kunita's theory on
stochastic flows of diffeomorphisms (see \cite{kunita} Corollary 4.6.6
for the first statement and equation (21) in the proof of Theorem 4.6.5
for the second).
\begin{aprop}
Let $a,b:[0,T]\times\R\rightarrow \R$ be functions bounded together with
their derivatives with respect to their second variable up to the order $3$. Then the stochastic
   differential equation
\begin{align}
   &d\chi^{\xi}_t=a(t,\chi^{\xi}_t)dW_t+b(t,\chi^{\xi}_t)dt,\;t\leq T,\;\chi^{\xi}_0=\xi\label{eds}\\
&(\mbox{resp. }d\zeta^{\nu}_t=a(T-t,\zeta^{\nu}_t)d\hat{W}_t+(a\partial_2a-b)(T-t,\zeta^{\nu}_t)dt,\;t\leq T,\;\zeta^{\nu}_0=\nu)\label{edsret}\end{align}
admits a solution such that for each $(t,\omega)\in[0,T]\times \Omega$,
the map ${\xi}\rightarrow \chi^{\xi}_t$ (resp. ${\nu}\rightarrow
\zeta^{\nu}_t$) is an increasing diffeomorphism of $\R$. The derivative
$\partial_\xi \chi^{\xi}_t$ solves the SDE
\begin{equation*}
   d\partial_{\xi} \chi^{\xi}_t=\partial_2a(t,\chi^{\xi}_t)\partial_\xi
\chi^{\xi}_tdW_t+\partial_2b(t,\chi^{\xi}_t)\partial_\xi \chi^{\xi}_tdt,\;\partial_\xi \chi^{\xi}_0=1.
\end{equation*}
Moreover
\begin{equation}
   d\P\;a.s.,\;\forall (t,\xi)\in[0,T]\times\R,\;\zeta_t^{\chi^{\xi}_T}=\chi^{\xi}_{T-t}.\label{inverse}
\end{equation}\label{propflot}\end{aprop}
\begin{adem}[of Proposition \ref{propflot2}]
Under our assumptions, the functions $x\sigma(t,x)$, $x\beta(t,x)$ and
   $x(\sigma\partial_2\eta-\beta)(t,x)$ vanish for $x=0$ and are Lipschitz continuous in $x$
   uniformly for $t\in[0,T]$. Therefore existence and trajectorial uniqueness follows from standard result
   concerning SDEs.\\
For
$[\bar{a},\bar{b}](t,\xi)=\left[\sigma,\beta-\frac{\sigma^2}{2}\right](t,e^\xi)$,
one has
\begin{equation}
   [\bar{a}\partial_2\bar{a}-\bar{b}](t,\xi)=\left[\sigma
e^\xi\partial_2\sigma
+\sigma^2-\beta-\frac{\sigma^2}{2}\right](t,e^\xi)=\left[\sigma\partial_2\eta-\beta-\frac{\sigma^2}{2}\right](t,e^\xi).\label{changfonc}
\end{equation}By
It\^o's formula, one easily deduces that for $\bar{\chi}$ and
$\bar{\zeta}$ solving 
\begin{align*}
   &d\bar{\chi}^{\xi}_t=\bar{a}(t,\bar{\chi}^{\xi}_t)dW_t+\bar{b}(t,\bar{\chi}^{\xi}_t)dt+dL_t,\;t\leq
T,\;\chi^{\xi}_0=\xi\\
&d\bar{\zeta}^\nu_t=\bar{a}(T-t,\bar{\zeta}^\nu_t)d\hat{W}_t+\left(\bar{a}\partial_2\bar{a}-\bar{b}\right)(T-t,\bar{\zeta}^\nu_t)+d\hat{L}_t,\;t\leq
T,\;\bar{\zeta}^\nu_0=\nu.
\end{align*}
$(e^{\bar{\chi}^{\log
    x}_t})_{t\in[0,T]}$ solves the SDE satisfied by $X^x$ and $(e^{\zeta^{\log z}_{t^+}})_{t\in[0,T]}$ where by
convention $\zeta^{\nu}_{T^+}=\zeta^{\nu}_{T}$ solves the SDE satisfied by $Z^z$
as soon as $e^{\zeta^{\log z}_{0^+}}=z$. Since the last equality holds as soon as
$L_{T^-}=L_T$ and therefore with probability $1$, by trajectorial
uniqueness, $$\P\left(Z^z_T=e^{\zeta^{\log z}_{T}}\mbox{ and }\forall
  t\in [0,T],\;(X^x_t,Z^z_t)=(e^{\bar{\chi}^{\log
    x}_t},e^{\bar{\zeta}^{\log
    z}_{t^+}})\right)=1.$$
With \eqref{changfonc}
one deduces \eqref{soledsback}.\par
For a fixed realization of $\mu$ or equivalently of the L\'evy process
$(L_t)_{t\in[0,T]}$, setting
$[a,b](t,\xi)=[\bar{a},\bar{b}](t,\xi+L_t)$, the
process $\chi^\xi_t=\bar{\chi}^\xi_t-L_t$ and
$\zeta^\nu_t=\bar{\zeta}^{\nu+L_T}_t-L_{T-t}$ respectively
solve \eqref{eds} and \eqref{edsret}. Since the functions $a$ and $b$
satisfy the hypotheses in Proposition \ref{propflot}, this result
implies that
$\xi\rightarrow\chi^\xi_T=\bar{\chi}^\xi_T-L_T$ and $\nu\rightarrow
\zeta^\nu_T=\bar{\zeta}^{\nu+L_T}_T$ are inverse increasing
diffeomorphisms of $\R$ and
$$\partial_\xi\bar{\chi}^\xi_T=\partial_\xi\chi^\xi_T=e^{\int_0^T\partial_2 {a}(t, {\chi}^{\xi}_t)dW_t+\int_0^T(\partial_2 {b}-\frac{1}{2}(\partial_2 {a})^2)(t, {\chi}^{\xi}_t)dt}=e^{\int_0^T\partial_2\bar{a}(t,\bar{\chi}^{\xi}_t)dW_t+\int_0^T(\partial_2\bar{b}-\frac{1}{2}(\partial_2\bar{a})^2)(t,\bar{\chi}^{\xi}_t)dt}.$$
Since $\xi=\zeta_T^{\chi^\xi_T}=\bar{\zeta}^{(\bar{\chi}^\xi_T-L_T)+L_T}_T=\bar{\zeta}^{\bar{\chi}^\xi_T}_T$, $\xi\rightarrow \bar{\chi}^\xi_T$ and $\nu\rightarrow
\bar{\zeta}^{\nu}_T$ and therefore $x\rightarrow X^x_T$ and
$z\rightarrow Z^z_T$ are inverse increasing diffeomorphisms. Equality
\eqref{ensemblegaux} follows easily. Moreover, on the almost sure event $\{\forall t\in[0,T],\;X^x_t=e^{\bar{\chi}^{\log
    x}_t}\}$,
\begin{align*}
   \partial_x
   X^x_T=\frac{1}{x}\partial_\xi \bar{\chi}^{\log x}_Te^{\bar{\chi}^{\log
    x}_T}=e^{\int_0^T(\bar{a}+\partial_2\bar{a})(t,\bar{\chi}^{\log x}_t)dW_t+\int_0^T[\bar{b}+\partial_2\bar{b}-\frac{(\partial_2\bar{a})^2}{2}](t,\bar{\chi}^{\log x}_t)dt+L_T}.
\end{align*}
One deduces \eqref{edsder} by remarking that
$\bar{a}+\partial_2\bar{a}(t,\xi)=(\sigma+e^\xi\partial_2\sigma)(t,e^\xi)=\partial_2\eta(t,e^\xi)$
and
\begin{align*}
   \left[\bar{b}+\partial_2\bar{b}-\frac{(\partial_2\bar{a})^2}{2}\right](t,\xi)&=\left[\beta-\frac{\sigma^2}{2}+e^\xi\partial_2\beta-e^\xi\sigma\partial_2\sigma-\frac{(e^\xi\partial_2\sigma)^2}{2}\right](t,e^\xi)
\\&=\left[\beta+e^\xi\partial_2\beta-\frac{(\sigma+e^\xi\partial_2\sigma)^2}{2}\right](t,e^\xi)=\left[\beta+e^\xi\partial_2\beta-\frac{(\partial_2\eta)^2}{2}\right](t,e^\xi).
\end{align*}
\end{adem}
\begin{adem}[of Proposition \ref{propflot3}]
When for fixed $z>0$, the function $\sigma$ does not vanish in a
neighbourhood of $(T,z)$, then for some
   $\varepsilon\in(0,T)$ the function $a$ does not vanish on
   $[T-\varepsilon,T]\times[\log z-2\varepsilon,\log
   z+2\varepsilon]$. Conditionally on $(L_t)_{t\in[0,T]}$, as soon
   as $L_{T^-}=L_T$, by
   Bouleau and Hirsch absolute continuity criterion (see Theorem 2.1.3
   p.162 \cite{bouhi}) , $\bar{\zeta}^{\log z}_T=\zeta^{\log z-L_T}_T$ has a density with respect to the
   Lebesgue measure and for all $\xi\in\R$, either
   $\P(\bar{\chi}^\xi_T\in[\log z-\varepsilon,\log z+\varepsilon])=\P(\chi^\xi_T\in[\log z-L_T-\varepsilon,\log z-L_T+\varepsilon])=0$ or the conditional law
   of $\bar{\chi}^\xi_T=\chi^\xi_T+L_T$ given $\chi^\xi_T\in[\log z-L_T-\varepsilon,\log z -L_T+\varepsilon]$ has a
   density. As $\P(L_{T^-}=L_T)=1$, one deduces \eqref{points}.\\
Let us now suppose that $\beta(t,x)=\gamma$. Then $X^x_t=xe^{\int_0^t\sigma(s,X^x_s)dW_s+\gamma
  t-\frac{1}{2}\int_0^t\sigma^2(s,X^x_s)ds+L_t}$ and
\begin{align*}
 {\mathbb E}(e^{-\gamma t}X^x_t)=x\E\left(e^{L_t}\E\left(e^{\int_0^t\sigma(s,X^x_s)dW_s-\frac{1}{2}\int_0^t\sigma^2(s,X^x_s)ds}\bigg|(L_s)_{s\in[0,T]}\right)\right).
\end{align*}
Since $\sigma$ is bounded, the conditional expectation in the
right-hand-side is equal to $1$ by Novikov's criterion (Proposition
1.15 p.308 \cite{reyor}) and the right-hand-side is equal to $x$ as
$\E(e^{L_t})=1$ (take $u=-i$ in \eqref{levykin}). Moreover, as the
function $\partial_2\eta(t,x)=(\sigma+x\partial_2\sigma)(t,x)$ is bounded, the expectation of $e^{-\gamma T}\partial_x X^x_T=e^{\int_0^T\partial_2\eta(s,X^x_s)dW_s-\frac{1}{2}\int_0^T(\partial_2\eta)^2(s,X^x_s)ds+L_T}$ is equal to $1$ by the
same argument.\\
Let $(h_n)_{n\geq 0}$ be a sequence of non-zero real numbers greater
   than $-x$ converging to
   zero. The random variables
   $\left((X^{x+h_n}_T-X^x_T)/h_n\right)_{n\geq 0}$ converge
   to $\partial_x X^x_T$ as $n\rightarrow +\infty$. Since they are
   non-negative and
$$\lim_{n\rightarrow
  +\infty}\E\left(\frac{X^{x+h_n}_T-X^x_T}{h_n}\right)=e^{\gamma T}=\E\left(\partial_x X^x_T\right),$$
they converge in $L^1$ to $\partial_x X^x_T$, which ensures uniform integrability.
\end{adem}
\end{document}